\documentclass[reqno, 11pt]{amsart}
\usepackage{amsmath,amssymb,supertabular}
\def\n{\noindent}
\def\hf{\hfil\break}

\begin{document}

\centerline {\bf A more intuitive definition of limit}
\smallskip

\centerline {Bogdan M. Baishanski}
\bigskip

\n {\bf Abstract}
To define -at the most elementary level- the right and the left limit of
the function $f$ at a real number $c$ (or the limit at $\infty$ or $-\infty$),
we propose these four axioms: if $f$ is ultimately monotone and bounded,
$\lim f$ exists; if $f =$ constant $= c$, then  $\lim f = c$; the sandwich
theorem holds; and, if  $\lim f$ is strictly less than  $\lim g$, then
ultimately  $f$ is strictly less than  $g$.  (In that elementary exposition,
only a few of results below would be included: at most the example, lemmas
1 and 2, and theorems 3',4,6 and 7).

This note is meant for a second, more advanced level, here it is shown
that the limit, i.e. a unique mapping satisfying the four axioms, exists.

There is also a third level, where this topic is considered in a more
abstract framework. We have given such an exposition in ``Axiomatic
definition of limit''.
 
\bigskip

\n {\bf Introduction}

The currently used $\epsilon-\delta$ definition of the limit appears to many
beginners to be too difficult. In the rearrangement of the introductory
material that we now propose, the epsilon-delta definition would be
dethroned and would get the status of a technical result (see thms 7 and 8
below). The new definition is obtained by combining two well-known
mathematical facts, which in current expositions appear as trivial
consequences of the $\epsilon-\delta$ definition. One of these facts (see (2)
below) becomes now the fundamental property of the limit.

In this note we consider real-valued functions of a real variable and define
lim $f(x)$ as $x$ tends to infinity. With obvious modifications one can phrase
definitions in the cases $x$ tends to $c+$, to $c-$, to $-\infty$, and similarly
define limit of a sequence of real numbers. Also, only minor modifications
are necessary to define in these cases the extended limit (that can take
as values real numbers, $\infty$ and  $-\infty$).

According to the $\epsilon-\delta$ definition, limit is a relation between
functions and numbers; according to our definition in this note, limit is
a mapping from a class of real functions to real numbers. So we need now
to define both the mapping and the maximal class on which that mapping can
be defined so that the defining property (2) of the limit is satisfied. We
do that in three steps: 

\n 1) we define limit of monotone bounded functions,
                    
\n  2) we define class of convergent functions (new
meaning is attached now to the term convergent- later it becomes evident
that it coincides with the traditional meaning),

\n 3) we extend the definition of the limit to all
convergent functions.

(It follows easily from theorems 7 and 8 that the class of convergent
functions cannot be further extended with preservation of property (2)).

Is our definition more accessible and more intuitive than the
epsilon-delta definition? Before answering that question, the reader
should consider possible graphical illustrations of definitions 1
and 2 below and compare them with graphical illustration of the
$\epsilon-\delta$ definition.

\n {\bf Results and proofs}

We start by considering functions which are ultimately monotone and
bounded,i.e. we consider the class $BM(\infty) = \{f|$  there exists $x _0$ such that $f$ is bounded and monotone on the interval $(x_0,\infty)\}$.

\n {\bf Def.1}.   A mapping
          $$L:   BM(\infty) \to  R$$
is said to be a limit on $BM(\infty)$ if the following two conditions are
satisfied:

\n (1) if $f(x) = c$   for every $x$, then   $L(f) = c$,

\n (2) if  $L(f) < L(g)$, then there exists a such that
               $$ f(x) < g(x)~{\rm  for}~  x > a.$$
Instead of  $L(f) = \lambda$   we shall frequently write  ``$\lim f(x) =\lambda$''   or
``$f(x)\to\lambda  ~{\rm  as}~ x\to \lambda$   '' . In case $f$ is increasing/decreasing we may write
$f(x)\nearrow \lambda$    as $x \to  \infty$ , resp. $f(x)\searrow \lambda $  
as $x\to \infty$.

\n {\bf Thm.1}.  There exists a limit $L$  on  $BM(\infty)$.

\n  {\bf Proof}. If $f$ is bounded and increasing for $x > a$, we set
$L(f) = \sup {f(x)| x > a}$, if f is bounded and decreasing for $x > a$, we set
$L(f) = \inf f(x)| x > a$. We need only to verify that the condition (2) is
satisfied. There are four cases to consider since  $f$ can be increasing or
decreasing,  and  similarly $g$. We shall consider here only the case when
$f$ is decreasing for $x > a'$ and $g$ is increasing for $x > a''$. (the other
three cases are even simpler).  Choose $\gamma$ such that $L(f)< \gamma < L(g)$. Since
$L(f) = \inf f(x)| x \ge a', L(g) = \sup g(x)| x\ge a''$, we obtain that there
exist $b'$ and $b"$ such that $f(b') <\gamma < g(b")$. Since $f$ is decreasing and $g$
is increasing we get that $f(x) < g(x)$ for $x > \max (b',b")$.

\n {\bf Thm.2.} There is only one limit $L$  on  $BM(\infty)$.

\n {\bf Proof.}  Suppose that there exist two limits, $L'$ and $L"$, and a function $f$
in $BM(\infty)$ such that $L'(f)$ is different from $L"(f)$, say  $L'(f) <
L"(f)$. Let $\gamma$ be a real number such that

\n  (3)                     $ L'(f) < \gamma < L"(f).$
Using $\gamma$ to denote also the constant function with value $\gamma$, we have by (1) that  $L'(\gamma) = L"(\gamma) = \gamma$ so from (3) it follows that $ L'(f) < L'(\gamma)$ and
$L"(\gamma) < L"(f)$, and therefore by (2) there exist $a'$ and $a"$ such that
$f(x) < \gamma$ for $x>a'$, and $f(x) >\gamma$ for $x>a"$. So, if $x > \max(a',a")$, we get $f(x) < f(x)$, a contradiction.

\n {\bf Thm. 3.}   If $f$ and $g$ belong to $BM(\infty)$ and if there exists $a$ such that
                    $$ f(x)\le g(x)~{\rm for}~ x>a, ~{\rm then}~
                          L(f) \le L(g).$$
                          
\n {\bf Proof.}  Assume not true. Then  $L(f) > L(g)$. By (2) it follows that there
exists $a'$ such that $f(x)> g(x)$ for $x>a'$. Taking an $x > max(a,a')$ we obtain
a contradiction.

\n {\bf Thm. 4.}   Let $N$ be a positive decreasing function on $[a,\infty)$. Then

\n (i) $L(N) = 0$
if and only if

\n (ii) for every positive integer n there exists $x_n > a$ such that
$N(x_n) < 1/n$.

\n {\bf Proof.}  (i)-$\to$(ii).   Since $L(N) < L(1/n)$ we get from (2) that there exists $c$
such that $N(x) < 1/n$  for $x>c$.

(ii)$\to$(i).  Since $N$ is decreasing on $(a,\infty)$ we have that $0< N(x)< 1/n$
for $x> x_n$.  By  Thm. 3 and by (1) we get  $0 \le L(N) \le 1/n$ for every $n$.
Thus $L(N) = 0$.

\n {\bf Example}. For every $K>0, c>o$ the function $K x^{-c} \searrow 0$ as $x\to \infty$.

\n {\bf Corollary 1.}  Let $f(x) = \lambda + N(x), g(x) = \lambda - N(x)$ .The following
three statements are equivalent: 

\n (i) $N(x) \searrow 0$, 

\n (ii) $f(x)\searrow \lambda$, 

\n (iii) $g(x) \nearrow \lambda$.

\n {\bf Corollary 2.}   Let $N'(x)\searrow 0, N"(x)\searrow 0$ as $x\to\infty$.

\n (i) if $N(x) = N'(x)+ N"(x)$, then  $N(x)\searrow 0$ as $x\to\infty$.

\n (ii)  if $N(x) = C N'(x)$, where $C$ is a constant, then  $N(x)\searrow 0$ as $x\to\infty$.

\n {\bf Definition 2} . We say that the function $f$ is convergent as $x$ tends to
$\infty$  if there exist a number $a$ and functions $m$ and $M$, both bounded
and monotone on the interval $[a, \infty)$ such that\hf
                        $f$ is defined on $[a, \infty)$\hf
                        $m(x) \le f(x) \le M(x)$, and\hf
                        $L(m) = L(M)$.\hf
We denote by $C(\infty$) the class of functions $f$ convergent as $x$ tends to
infinity.

\n {\bf Thm. 5}.  The limit $L$  on $BM(\infty)$ can be extended to the class $C(\infty)$,
in such a way that the condition (2) remains satisfied.

\n {\bf Proof.}   If $f$ belongs to $C(\infty)$, set  $L(f) = L(m) = L(M)$.
First we need to show that this is a definition of $L(f)$, i.e. if
$m'(x) \le f(x) \le M'(x)$ for $x>a' , m"(x) \le f(x) \le M"(x)$ for $x>a",
L(m') = L(M') = L'$, and $L(m") = L(M") = L"$,  then  $L' = L"$.
We derive from the assumptions above that $m'(x) \le M"(x)$ for $x>\max(a',a")$.
Thus by Thm. 3  $L' = L(m') \le L(M") = L"$. Similarly $L" \le L'$. So $L'= L"$.

It remains to be shown that the condition (2) is satisfied on $C(\infty)$.
Let this time  $m',M', m", M"$  belong to $BM(\infty)$ and satisfy
$m'(x) \le f(x) \le M'(x)$ for $x>a' , m"(x) \le g(x) \le M"(x)$ for $x>a",
L(m') = L(M') = L(f)$, and $L(m") = L(M") = L(g)$. Since
$L(M') = L(f) < L(g) = L(m")$ we obtain from (2) and Theorem 1 that there exists $a$ such that
$M'(x)< m"(x)$  for $x>a$. Therefore  $f(x) \le M'(x) < m"(x) \le g(x)$ for
$x>max(a,a',a")$.

\n {\bf Thm. 2'}.  The limit $L$  on $C(\infty)$ is unique.

Proof is identical to the proof of Thm. 2.

\n {\bf Thm. 3'}.   If there exists $a$ such that
                     $$ f(x)\le g(x)~{\rm  for}~ x>a,$$
and if $f$ and $g$ converge as $x$ tends to $\infty$, then
                      $$L(f) \le L(g).$$

\n {\bf Proof}.   Same as of Thm. 3.

\n {\bf Lemma 1}.

\n a) Let $|f(x)-\lambda | < N(x)$ for $x > a$. If $N(x) \searrow 0$ as $x \to
\infty$, then $f(x) \to \lambda$ as $x\to\infty$.

\n b)  Let $f(x) = \lambda + z(x)$. Then the following statements are
equivalent:

\n (i) $f(x)\to \lambda$ as $x\to \infty$ and

\n (ii) $z(x)\to 0$ as $x\to\infty$.

\n {\bf Lemma 2}.  (i)  If  $z'(x)\to 0, z"(x)\to 0$ as $x\to\infty$,  and $z(x)=z'(x)+z"(x)$,
then $z(x)\to 0$ as $x\to \infty$.

\n (ii)  If  $z(x)\to 0$ as $x\to\infty$, and $w(x)=b(x)z(x)$, where the
function $b$ is bounded on some interval  $[a,\infty)$, then
$w(x)\to 0$ as $x\to\infty$.

\n {\bf Thm. 6}. Let $f(x)\to \alpha, g(x)\to\beta$ as $x\to\infty$.  If $s = f + g$,
$p = f g$  and $q = 1/g$, then, as $x\to\infty,  s(x)\to\alpha + \beta,
p(x)\to \alpha  \beta$, and if $\beta$ is different from $0, q(x) \to 1/\beta$.

\n {\bf Proof}.   Follows easily from Lemmas 1 and 2.

\n {\bf Thm. 7}.   Let  $f$ converge as $x$ tends to $\infty$ and let $L(f) = \lambda$.

\n (i)  If  $\alpha < \lambda < \beta$, there exists $a$ such that
                   $\alpha < f(x) < \beta$ for $x > a$, and
                   
\n (ii)  If $\epsilon$ is a positive real number, there exists
$X=X(\epsilon)$ such that

                   $$|f(x)-\lambda| < \epsilon~{\rm if}~ x > X.$$

\n {\bf Proof}.  (i)   We shall use the same letter to denote a real number  and the
constant function which has that number as its unique value, and so we
write  $L(\alpha)= \alpha, L(\beta)= \beta$. The assumption is
               $$L(\alpha) < L(f) < L(\beta).$$
Applying twice (2) we get
            $$\alpha < f(x)~{\rm for}~ x > a', f(x) < \beta~{\rm for}~ x > a",$$
so that conclusion of (i) holds with $a=max(a',a")$.

\n (ii) Apply (i) with  $\alpha = \lambda - \epsilon, \beta = \lambda + \epsilon$.

\n {\bf Thm. 8}   If for every positive real number $\epsilon$ there exists $X =
X(\epsilon)$ such that

\n (4)  $|f (x)-\lambda | < \epsilon$ if $x > X$,
then (in the sense of Definition 2) $f$ converges as $x$ tends to infinity and
$L(f) = \lambda$.

\n {\bf Proof}.   We need to show that conditions of definition 2 are satisfied.
Let  $M(x) = \sup {f(t)| t \ge x}, m(x) = inf {f(t)| t \ge x}$.  It is easy to
verify that the functions $M$ and $m$ belong to $BM(\infty)$ and that
$m(x) \le f(x) \le M(x)$, say for $ x > X(1)$. We need only to verify that
$L(m) = L(M) =\lambda$. We shall show that by proving $L(M) =\lambda$,
analogously one can prove $L(m) = \lambda$.

Given a positive real number $\epsilon$, we derive from (4) that
$\lambda - \epsilon < f(t) <  \lambda + \epsilon$ if $t > X(\epsilon)$. Therefore
$\lambda - \epsilon \le M(x) \le  \lambda + \epsilon$ if $x > X(\epsilon)$. Applying
Thm. 3 we get that, for every positive $\epsilon~\lambda - \epsilon \le L(M) \le  \lambda + \epsilon$.   
That implies $L(M) = \lambda$.

\bigskip
 
\n Bogdan Baishanski\hf
Department of Mathematics\hf
Ohio State University\hf
231 West 18 Avenue\hf
Columbus, OHIO 43210
\hf
bogdan@math.ohio-state.edu
\end{document}